\documentclass[12pt]{article}
\usepackage{amsfonts}
\usepackage{amsthm}
\usepackage{latexsym}
\usepackage{latexsym}
\usepackage{amsthm}
\usepackage{graphicx}
\makeatletter
\renewcommand{\@seccntformat}[1]
{\csname the#1\endcsname.\enspace} \makeatother
\renewcommand{\baselinestretch}{1.75}
\begin{document}
\newtheorem{theorem}{Theorem}
\newtheorem{lemma}{Lemma}
\newtheorem{remark}{Remark}
\newtheorem{remarks}{Remarks}
\newtheorem{corollary}{Corollary}
\newtheorem{example}{Example}
\newtheorem{definition}{Definition}
\newtheorem{assumption}{Assumption}
\thispagestyle{empty}

\begin{center}
   {\bf On Bayesian credible sets in restricted parameter space problems and lower bounds for frequentist coverage \footnote{\today} }
\end{center}

\begin{center}
{\sc \' Eric Marchand$^{a}$, William E. Strawderman$^{b}$} \\

{\it a  Universit\'e de
    Sherbrooke, D\'epartement de math\'ematiques, Sherbrooke Qc,
    CANADA, J1K 2R1 (e-mail: eric.marchand@usherbrooke.ca) } \\

{\it b  Rutgers University, Department of Statistics and Biostatistics, 501 Hill
Center, Busch Campus, Piscataway, N.J., USA, 08855 (e-mail:
straw@stat.rutgers.edu) }
\end{center}

\begin{center}
{\sc Summary}
\end{center}
\small 
For estimating a lower bounded parametric function in the framework of Marchand and Strawderman (2006), we provide through a unified approach a class of Bayesian confidence intervals with credibility $1-\alpha$ and frequentist coverage probability bounded below by $\frac{1-\alpha}{1+\alpha}$.  In cases where the underlying pivotal distribution is symmetric, the findings represent extensions with respect to the specification of the credible set achieved through the choice of a {\it spending function}, and include Marchand and Strawderman's HPD procedure result.  For non-symmetric cases, the determination of a such a class of Bayesian credible sets fills a gap in the literature and includes an ``equal-tails'' modification of the HPD procedure.   Several examples are presented demonstrating wide applicability. 


\normalsize
\section{Introduction}
Bayesian credible sets are not designed (e.g., Robert, 2011) and are far from guaranteed (Fraser, 2011) to have
satisfactory, exact or precise frequentist coverage but it is nevertheless of interest to investigate (Wasserman, 2011) to what
extent there is convergence or divergence in various situations.
A historically resonating example where there is exact convergence arises for estimating the mean of a $N(\mu, \sigma^2)$ distribution, and where the use of the non-informative prior leads to a $(1-\alpha) \times 100\%$ HPD credible set (i.e. the $z$ or $t$ confidence interval) with exact frequentist coverage.  This, however, is very much the exception.  Even, in the simple presence of a lower bound on the mean parameter $\mu$ (e.g., Mandelkern, 2002), with the prior taken to be the truncation of the non-informative prior onto the restricted parameter space, the frequentist coverage of the $(1-\alpha) \times 100\%$ HPD credible set fluctuates from its credibility (or nominal coverage) $1-\alpha$.  However, the HPD procedure does not fare poorly as a frequentist procedure for large $1-\alpha$ as witnessed by the lower bound $\frac{1-\alpha}{1+\alpha}$ on its frequentist coverage due to Roe and Woodroofe (2000, known $\sigma^2$) and Zhang and Woodroofe (2003, unknown $\sigma^2$), 
as well as the better lower bound $1 - \frac{3\alpha}{2}$ (for $\alpha < 1/3$, known $\sigma^2$) obtained by Marchand et al. (2008).  

In a generalization of the above, Marchand and Strawderman (MS 2006) introduced a unified framework for which the 
$(1-\alpha) \times 100\%$ HPD credible set of a lower bounded parametric function has
frequentist coverage greater than $\frac{1-\alpha}{1+\alpha}$ for all values lying in the restricted parameter space.  
This framework, as well as its various applications, will be revisited in Sections 2 and 5, but let us consider for sake of illustration the basic examples:
(i) $X \sim f_0(x-\theta)$ with known $f_0$, $\theta \geq 0$; and (ii) $X \sim \hbox{Gamma}(\alpha, \theta)$ with $\theta \geq 1$ with known $\alpha$.
For location family densities as in (i) with $f_0$ unimodal and symmetric, Marchand and Strawderman's results apply for 
the flat prior on the truncated parameter space $[0,\infty)$ and the corresponding $(1-\alpha) \times 100\%$ HPD credible set, with the guarantee that the actual frequentist coverage is bounded below by $\frac{1-\alpha}{1+\alpha}$ for all $\theta \geq 0$.  However, if $f_0$ is not symmetric, such a result does not hold in 
general (MS 2006, Example 1).  The same is true for a vast number of so-called non-symmetric situations arising in Marchand and Strawderman's framework, including 
the Gamma models in (ii) where the prior is given by $\frac{1}{\theta} \mathbb{I}_{[1,\infty)}(\theta)$, that is the truncation on $[1,\infty)$ of the usual non-informative prior $\frac{1}{\theta} \mathbb{I}_{(0,\infty)}(\theta)$.   It is true that the bound holds for certain specific classes of $f_0$'s (MS 2006, Theorem 2, a), and it is also the case that numerical evaluations of a theoretical and unexplicit lower bound for frequentist coverage provides further evidence for satisfactory coverage
for a specific Gamma model in (ii) (MS 2006, Example 2).    Nevertheless, a clear analytical result or lower bound for
frequentist coverage in such non-symmetric cases is lacking, and it our motivation here to try to fill this gap. 
 
For a large variety of situations with a lower bounded parametric constraint, we obtain here a class of Bayesian $(1-\alpha) \times 100\%$ 
credible sets which provide minimal frequentist probability coverage exceeding $\frac{1-\alpha}{1+\alpha}$.  These Bayesian confidence intervals include
an ``equal-tails'' modification, or approximation, of the HPD credible set, which also coincides with the latter in situations of underlying symmetry.  
Our findings are achieved by introducing and exploiting a {\it spending} function interpretation of Bayesian confidence intervals, and lead to a class of procedurees 
(rather than a single one) which share the above lower bound for frequentist coverage.   The rest of the paper is organized as follows.  Preliminary results, definitions and model assumptions, including those related to the spending function associated with a Bayesian credible interval,  are presented in Section 2, while Bayesian credible interval representations are outlined in Section 3.   The main findings concerning frequentist coverage appear in Section 4 and various 
applications are presented and commented on in Section 5. 

\section{Definitions and preliminary results }

\subsection{Assumptions, invariance, pivot, prior, and implications}

As in basic examples (i) and (ii), we consider model densities $f(x;\theta)$; $x \in \cal{X}$, $\theta \in \Theta \subset \mathbb{R}^p$; for an
observable $X$, and we are concerned with interval estimation of a parametric function $\tau(\theta)$
($\mathbb{R}^p \to \mathbb{R}$) with the additional constraint $\tau(\theta) \geq 0$.    We assume there exists a pivot
of the form $T(X,\theta)=\frac{a_1(X)-\tau(\theta)}{a_2(X)}$; $a_2(\cdot)>0$;  such that $-T(X,\theta)$ has cdf $G$ and Lebesgue density $g_0$.  
This pivot assumption means that the frequentist or
conditional distribution of $T(X,\theta)$, or $-T(X,\theta)$, given $\theta$ does not depend on $\theta$, $\theta \in \mathbb{R}^{p}$.
We can thus set $G$ as the common cdf of $-T(X,\theta)$.  In the basic location-family example (i) with $X \sim f_0(x-\theta)(=g_0(\theta-x), \hbox{ say})$, the above is illustrated by the fact that $-T(X,\theta)=\theta-X$ is a pivot with cdf $G$ and pdf $g_0$.  In the Gamma example, or more generally scale families with $X \sim \frac{1}{\theta} \, f_1(\frac{x}{\theta}), \theta \geq 1$, a corresponding $T(X,\theta)$ pivot is obtained with $a_1(X)=\log(X)$, $a_2(X)=1$, $\tau(\theta)=\log(\theta)$.

We further assume that the unrestricted decision problem is invariant under a group $\cal{G}$ of transformations and that the pivot satisfies the invariance requirement
$T(x,\theta)=T(gx,\bar{g}\theta)$, for all $x \in \cal{X}$,
$\theta \in \Theta$, $g \in \cal{G}$, $\bar{g} \in \bar{\cal{G}}$, with
$\cal{X}$, $\Theta$, $G$, and $\bar{G}$ being isomorphic.  For instance, in basic example (i), the invariance is achieved with the additive group $G$ on $\mathbb{R}^p$
and since $T(x,\theta)=x-\theta = (x+g)-(\theta+g)=T(gx,\bar{g}\theta)$ for all group elements $g$. 

Collecting the above assumptions, we have for further reference.

\begin{assumption}
\label{assumption}
We have a model density $f(x;\theta)$; $x \in \cal{X}$, $\theta \in \Theta$; for an
observable $X$, with both $X$ and $\theta$ being vectors, and we
seek to estimate a parametric function $\tau(\theta)$ $(\mathbb{R}^p \to \mathbb{R})$ with the
constraint $\tau(\theta) \geq 0$.    We assume there exists a pivot
$T(X,\theta)=\frac{a_1(X)-\tau(\theta)}{a_2(X)}$; $a_2(\cdot)>0$;
such that $-T(X,\theta)$ has cdf $G$ and Lebesgue density $g_0$. 
 We further assume that the decision problem is invariant under a group $\cal{G}$ of transformations and that the pivot satisfies the invariance requirement
$T(x,\theta)=T(gx,\bar{g}\theta)$, for all $x \in \cal{X}$,
$\theta \in \Theta$, $g \in \cal{G}$, $\bar{g} \in \bar{\cal{G}}$, with
$\cal{X}$, $\Theta$, $G$, and $\bar{G}$ being isomorphic.   
\end{assumption}

We consider prior measures $\pi_H$ and $\pi_0$, where $\pi_0(\theta) =
\pi_H(\theta) I_{[0, \infty)}(\tau(\theta))$, and $\pi_H$ is the Haar
right invariant measure which satisfies the property $\pi_H( A \, \bar{g}) = \pi_H(A)$ for every measurable subset $A$ of $\Theta$, and for every $g \in G$.  The right Haar measure $\pi_H$ exists and is unique up to a multiplicative constant for locally compact groups such as location, scale, and location-scale.  For the basic location and the Gamma model (or scale model) examples of the Introduction,  right Haar invariant measures are given by $\pi_H(\theta)=1$ and $\pi_H(\theta)=\frac{1}{\theta}$ respectively.  For a sample from a location-scale family with 
$X_i  \sim^{\hbox{ind.}} \frac{1}{\theta_2} f_2(\frac{x_i-\theta_1}{\theta_2})\;, i=1, \ldots,n$,
the common non-informative prior $\pi(\theta)= \frac{1}{\theta_2}$ is right Haar invariant.  We refer to Berger (1985) or Eaton (1989) for detailed treatments of invariance and Haar invariant measures.
 
A key feature relative to Assumption 1 and the choice of the right Haar invariant measure is that the frequentist distribution
of $T(X,\theta)$; which is free of $\theta$ by virtue of the pivot assumption for  $T(X,\theta)$; coincides with the posterior distribution of $T(x,\theta)$ under $\pi_H$ for any given $x$, i.e., 
\begin{equation}
\label{t=t} 
T(x,\theta) | x \,\; =^d \; T(X,\theta)|\theta\,,   \hbox{ for all } x, \theta.
\end{equation}
We will pursue, after the next Lemma, by illustrating the above and drawing implications of immediate interest.  For sake of completeness, we reproduce here a key lemma from MS(2006) justifying (\ref{t=t}) and we refer to their work for further  details.

\begin{lemma} (MS, 2006, Corollary 1)
\label{haar2} Suppose $\cal{X}$, $\Theta$, $G$, and $\bar{G}$ are
all isomorphic, and that $T(X,\theta)$ is a function for which
$T(x,\theta)=T(gx,\bar{g}\theta)$, for all $x \in \cal{X}$,
$\theta \in \Theta$, $g \in G$, $\bar{g} \in \bar{G}$.  Then
condition (\ref{t=t}) holds, that is $P_{\theta}[T(X,\theta) \in
B] = P^{\pi_H(\theta|x)}[T(X,\theta) \in B]$ for each measurable
set $B$. 
\end{lemma}

 Now, for the basic unrestricted location family example with the flat prior $\pi_H(\theta)=1$, which is Haar right invariant, observe that the posterior density of $\theta$ is given by
$$\pi_H(\theta|x) = \frac{f_0(x-\theta)}{\int_{\theta} f_0(x-\theta) \, d\theta} = f_0(x-\theta) = g_0(\theta-x)\,,  $$ 
so that the posterior density of $-T(x,\theta)=\theta-x$ associated with $\pi_H$ is given by $g_0$ as well.
This correspondence for basic example (i) illustrates property (\ref{t=t}) which is, of course, more general under Assumption 1.

In general, observe that the posterior cdf under $\pi_H$ for $\tau(\theta)$ is available from the fact that $-T(X,\theta)=\frac{\tau(\theta)-a_1(X)}{a_2(X)} \sim G$ yielding 
\begin{equation}
\label{cdfunderpi}
P_{\pi_H}(\tau(\theta) \leq y |x) = G(\frac{y-a_1(x)}{a_2(x)})\,.
\end{equation}
Now, under the truncation $\pi_0$ of $\pi_H$, the above correspondence  between the frequentist and posterior distributions of $-T(X,\theta)$ does not hold, and the posterior cdf under $\pi_0$ of $\tau(\theta)$ differs.  However, we can still express the posterior distribution of $\tau(\theta)$ under $\pi_0$ in terms of $\pi_H$ and $G$.  Indeed, with $\pi_0(\theta)=\pi_H(\theta) \, \mathbb{I}_{[0,\infty)}(\tau(\theta))$ and $\frac{\tau(\theta)-a_1(x)}{a_2(x)}|x \sim G$ under $\pi_H$, we have for
a measurable set $A \subset \Theta_0=\{\theta \in \Theta : \tau(\theta) \geq 0 \}$, and for any $x$:
\begin{eqnarray*}
P_{\pi_0}(\theta \in A |x) &=& \int_A \, \pi_0(\theta|x) \, d\theta \\
\, &=& \frac{\int_A \, \pi_0(\theta) \, f(x|\theta) \, d\theta }{\int_{\Theta_0} \, \pi_0(\theta) \, f(x|\theta) \, \, d\theta } \\
\, &=& \frac{\int_A \, \pi_H(\theta) \, f(x|\theta) \, d\theta }{\int_{\Theta_0} \, \pi_H(\theta) \, f(x|\theta) \, \, d\theta }\\
\, &=&  \frac{P_{\pi_H}(\theta \in A|x)}{P_{\pi_H}(\theta \in \Theta_0|x)}\,. 
\end{eqnarray*}
In terms of the posterior survival function of $\tau(\theta)$ under $\pi_0$, the above yields along with (\ref{cdfunderpi}), for $y \geq 0$,
\begin{equation}
\label{survival} P_{\pi_0}(\tau(\theta) \geq y |x) \,=\, \frac{P_{\pi}(\tau(\theta) \geq y|x)}{P_{\pi}(\tau(\theta) \geq 0|x)}\,
= \, \frac{1-G(\frac{y-a_1(x)}{a_2(x)})}{1-G(\frac{-a_1(x)}{a_2(x)})} \,. 
\end{equation}

We will make use, in Section 3, of the above in setting and describing the bounds of Bayesian credible sets for $\tau(\theta)$ under $\pi_0$.

\subsection{The spending function associated with a Bayesian credible set}

With the objective of constructing a $(1-\alpha) \times 100 \%$ Bayesian credible set or region, the determination of a posterior distribution
for $\tau(\theta)$ supported on $[0,\infty)$ leaves open many choices and various different approaches (e.g., Berger, 1985, section 4.3.2).  The HPD
credible set is one such region chosen to minimize volume and leading to intervals for unimodal posterior densities.   In our set-up, $(1-\alpha) \times 100 \%$
Bayesian credible intervals are, more generally, of the form $[l(x), u(x)]$, $x \in \mathbb{R}$, where $P(l(x) \leq \tau(\theta) \leq u(x)|x)=1-\alpha$.
An alternative (and equivalent) way to set or view the bounds $l(x)$ and $u(x)$, for a given $x$, is to focus on the complementary set $[0,l(x)) \cup (u(x), \infty)$ 
and to allocate (or ``spend'') probabilities $\alpha-\alpha(x)$ and $\alpha(x)$ respectively on its two disjoint parts, with $\alpha(x) \in [0,\alpha]$.  
It is clear (when the posterior density is absolutely continuous) that the choice $\alpha(x)$ leads to a unique choice of $[l(x), u(x)]$, and vice-versa.  Since we are interested in the frequentist properties of such Bayesian credible intervals, we will represent this allocation as a {\it spending} function.  Moreover, our findings guaranteeing minimal frequentist coverage of at least $\frac{1-\alpha}{1+\alpha}$ for a class of Bayesian credible sets will be conveniently expressed as conditions on the corresponding spending function.

\begin{definition}
\label{spendingdefinition}
For a given prior $\pi$ for $\theta $ and a credibility coefficient $1-\alpha$, a spending function $\alpha(\cdot): \mathbb{R}^p \to [0,\alpha]$ is a function such that, for all $x$, $P_{\pi}(\tau(\theta) \geq u(x)|x)=\alpha(x)$, $P_{\pi}(\tau(\theta) \leq l(x)|x)=\alpha-\alpha(x)$, and $[l(x),u(x)]$ is a 
$(1-\alpha) \times 100 \%$ Bayesian credible interval for $\tau(\theta)$.
\end{definition}
 
For example, a lower-tailed credible interval for a given $x$ corresponds to the selection $\alpha(x)=\alpha$, an upper tailed credible interval corresponds to $\alpha(x)=0$, and an equal tailed (based on the posterior $\pi$) corresponds to $\alpha(x)=\alpha/2$.

\subsection{Checklist}
 To facilitate the further presentation of
the results, here is a list of definitions and notations used.

{\it Cheklist} \baselineskip=0.5\normalbaselineskip

\begin{itemize}
\item  $1- \alpha$: credibility or posterior coverage or nominal
frequentist coverage ($\alpha \in (0,1)$)
\item  $T(X,\theta)=\frac{a_1(X)-\tau(\theta)}{a_2(X)}$:  pivot
\item $\pi_H$: unrestricted prior density chosen as the right Haar invariant measure
\item $\pi_0$ : prior density given by the truncation of $\pi_H$ onto the restricted parameter space 
\item  $G$: cumulative
distribution function (cdf) of $-T(X,\theta)|x$ and of $-T(X,\theta)|\theta$ under $\pi_H$ (which coincide for all $x,\theta$)
 \item $g_0= G'$: probability
density function (pdf) of $-T(X,\theta)$ 
\item $G^{-1}$: inverse cdf 
\item  $\alpha(\cdot)$: spending function
\item
$I_{\pi_0,\alpha(\cdot)}(X)=[l(X),u(X)]$: Bayesian credible set of
credibility $1-\alpha$ associated with the prior $\pi_0$ and the spending function
$\alpha(\cdot)$  \item
$C(\theta)$: the frequentist coverage at $\theta$ of the
confidence interval $I_{\pi_0,\alpha(\cdot)}(X)$ given by $C(\theta) =
P_{\theta}(I_{\pi_0,\alpha(\cdot)}(X) \ni \tau(\theta))$ 
\item  $y_0=-G^{-1}(\frac{\alpha}{1+\alpha})$
\item  $t(x)= \frac{a_1(x)}{a_2(x)} $
\item  $\Delta_0(x)=(1-\alpha)(1-G(-t(x)))$
\end{itemize}

\baselineskip=1.0\normalbaselineskip
\section{Bayesian credible intervals: representations and properties}

In this section, we expand upon  two different, yet equivalent, and instructive approaches to
constructing a credible set for $\tau(\theta)$ associated with prior $\pi_0$.  These are: (A) the spending function approach, and
(B) the approach based on the quantiles of the pivot.

\begin{enumerate}
\item[ {\bf (A)}]  (Spending function approach) \\
As seen above, a $(1-\alpha) \times 100\%$ credible interval for $\tau(\theta)$
associated with prior $\pi_0$ can be generated by a {\it spending}
function $\alpha(\cdot): \mathbb{R}^p \to [0,\alpha]$, such that
$I_{\pi_0, \alpha(\cdot)}(X)=[l(X), u(X)]$ with $P_{\pi_0}(\tau(\theta)
\geq u(x)|x) = \alpha(x)$. 
More precisely, we have the following under Assumption \ref{assumption}.

\begin{lemma}
\label{atol} For a given spending function $\alpha(\cdot)$, the bounds of 
$I_{\pi_0, \alpha(\cdot)}(x)$ are given by: $l_{\alpha(\cdot)}(x) = a_1(x) + a_2(x) G^{-1}
\{G(-t(x)) + (\alpha - \alpha(x)) (1-
G(-t(x) ) \}$ and $u_{\alpha(\cdot)}(x) = a_1(x)
+ a_2(x) G^{-1} \{ 1 - \alpha(x) (1- G(-t(x) ) \}
$, with $t(x)=\frac{a_1(x)}{a_2(x)}$.
\end{lemma}
{\bf Proof.}  With the survival function $ P_{\pi_0}(\tau(\theta) \geq y|\,x) = 
\frac{1-\,G(\frac{y-a_1(x)}{a_2(x)})}{1-G(-t(x))}\,$, as given in (\ref{survival}),
we obtain for $\beta \in (0,1), y >0$,
$P_{\pi_0}(\tau(\theta) \geq y|\,x) \,=\, \beta \; \Leftrightarrow
y = a_1(x) + a_2(x) G^{-1} \{1-\beta + \beta
G(-t(x))\}$, and the result follows with the
choices $\beta=\alpha(x)$ and $\beta=1-(\alpha - \alpha(x))$ for
$u(x)$ and $l(x)$ respectively.  \qed 

\begin{example}
\label{symmetric} The HPD procedures studied by MS (2006) for symmetric about $0$ and unimodal $g_0$ are given by the
bounds $l(x) = \max\{0, a_1(x) + a_2(x) G^{-1}(\frac{1-(1-\alpha)
G(t(x))}{2}) \}$ and $u(x) = a_1(x) + a_2(x)
\min\{G^{-1} (1-\alpha G(t(x))),
G^{-1}(\frac{1+(1-\alpha) G(t(x))}{2}) \}\,.$
With these given bounds, one may verify directly from
(\ref{survival}) that the corresponding spending function is
equal to
\begin{equation}
\label{spendinghpd} \min\{\alpha, \frac{\alpha}{2} +
\frac{G(-t(x)}{2(1- G(-t(x)))}
\},
\end{equation}
with $\alpha(x)=\alpha$ if and only if $t(x) \leq
-G^{-1}(\frac{\alpha}{1+\alpha})= G^{-1}(\frac{1}{1+\alpha})$
since $g_0$ is symmetric about $0$.  Conversely, applying
Lemma \ref{atol} with the spending function choice $\alpha(\cdot)$
in (\ref{spendinghpd}) leads to the HPD procedure above (using the
equality of $G(\cdot)$ and $1-G(-\cdot)$ for symmetric about $0$ $g_0$'s).
\end{example}

\item[ {\bf(B)}] (Approach based on quantiles of the pivot) \\
Alternatively, a second approach for cases where $l(x) >0$ begins
with choices $\gamma_1$ and $\gamma_2$, which will be made for each $x$, such that $G(\gamma_2) -
G(-\gamma_1) = \Delta$, for a given $\Delta \in (0,1)$.  Since,
for any $x$, we require $1-\alpha = P_{\pi_0}(l(x) \leq
\tau(\theta) \leq u(x)|x)$, we must have by (\ref{survival}):
$$G(\frac{u(x)-a_1(x)}{a_2(x)}) - G(\frac{l(x)-a_1(x)}{a_2(x)}) =
(1-\alpha)(1-G(-t(x))), $$ and this can be
achieved with choices $-\gamma_1$ and $\gamma_2$ above for
$\Delta= \Delta_0(x)=(1-\alpha)(1-G(-t(x)))$ yielding
$\frac{u(x)-a_1(x)}{a_2(x)} =
\gamma_2(\Delta_0(x))$ and
$\frac{l(x)-a_1(x)}{a_2(x)} =
-\gamma_1(\Delta_0(x))$, in other
words
\begin{equation}
\label{pivotbounds}
 l(x) = a_1(x) -a_2(x)\,\gamma_1
(\Delta_0(x)), \; \hbox{and } u(x) = a_1(x) +
a_2(x) \, \gamma_2(\Delta_0(x))\,,
\end{equation}
whenever $l(x) > 0$.  In view of the lower bound restriction on $\tau(\theta)$ (i.e., $\tau(\theta) \geq 0)$, and the corresponding
requirement that $l(X) \geq 0$,  observe that
not all choices of $-\gamma_1$ (and hence of $\gamma_2$) are feasible in (\ref{pivotbounds}) and that we must have
$$  -\gamma_1(\Delta_0(x)) \geq -\frac{a_1(x)}{a_2(x)}\,.$$
\end{enumerate}
\begin{example}
\label{equaltails}  With the above construction in (\ref{pivotbounds}), an {\it equal-tails} choice of $-\gamma_1$ and $\gamma_2$, that is $-\gamma_1(\Delta)=
G^{-1}(\frac{1-\Delta}{2})$ and $\gamma_2(\Delta)=
G^{-1}(\frac{1+\Delta}{2})$, leads to the credible interval bounds
\begin{equation}
\label{et} 
 l(x) = a_1(x) +a_2(x) \, G^{-1}(\frac{1-\Delta_0(x)}{2}), \; \hbox{and } u(x) = a_1(x) +
a_2(x)  \, G^{-1}(\frac{1+\Delta_0(x)}{2}),
\end{equation}
when $l(x) > 0$.  These above bounds coincide with those of the HPD procedure (when $l(x) > 0$) in
the symmetric case of Example \ref{symmetric}, as well as the
spending function given in (\ref{spendinghpd}) as can be verified
directly from (\ref{survival}). \\
NOTE:  We wish to emphasize that the terminology ``equal tails''
does not mean $\alpha(x)=\alpha/2$ (i.e., equal tails under the posterior distribution), but rather refers to the choice of (equal tails) 
quantiles $-\gamma_1$ and $\gamma_2$ under $G$.
\end{example}

The next section's lower bound of $\frac{1-\alpha}{1+\alpha}$ on frequentist coverage applies to a class of Bayesian credible intervals.
This class will include an {\it equal-tails} credible interval $I_{\pi_0, \alpha_{eqt}(\cdot)}$ which relates to both approaches presented in this section.  On one hand, it borrows the bounds (and hence the spending function) of the HPD procedure for symmetric about $0$ unimodal densities and, on the other hand, it is defined through the above equal-tailed choice (whenever $l(x)>0$).

\begin{definition}  
\label{etdef}
In the context of Assumption 1, the $G$-equal-tails credible interval $I_{\pi_0, \alpha_{eqt}(\cdot)}(X)$ is given by the bounds $l(x) = \max\{0, a_1(x) + a_2(x) G^{-1}(\frac{1-(1-\alpha) G(t(x))}{2}) \}$ and $u(x) = a_1(x) + a_2(x)
\min\{G^{-1} (1-\alpha G(t(x))),
G^{-1}(\frac{1+(1-\alpha) G(t(x))}{2}) \}\,.$  Equivalently, $I_{\pi_0, \alpha_{eqt}(\cdot)}(X)$ is given by the spending function
\begin{equation}
\label{alphaeqt} \alpha_{eqt}(x)= \min\{\alpha, \frac{\alpha}{2} +
\frac{G(-t(x)}{2(1- G(-t(x)))}
\}.
\end{equation}

\end{definition}

\section{Frequentist coverage properties}

We study here the frequentist coverage properties, under Assumption 1, of a class of Bayesian credible intervals which includes the equal-tails
credible interval $I_{\pi_0, \alpha_{eqt}(\cdot)}(X)$.  This procedure, as well as Example \ref{symmetric}'s HPD procedure for symmetric $g_0$, produces
estimates of the form $[0,u(x)]$ if and only if $t(x) \leq y_0$, where $y_0=-G^{-1}(\frac{\alpha}{1+\alpha})$ (and
$t(x)=\frac{a_1(x)}{a_2(x)}$ as above).  We thus focus on a class of credible intervals with the same behaviour.  Said otherwise in terms of the spending
function, we impose the choice $\alpha(x)=\alpha$ whenever $t(x) \leq y_0$.  We hence seek conditions on $\alpha(x)$, for those $x$'s such that
$t(x) \geq y_0$, for which minimal frequentist coverage is bounded below by $\frac{1-\alpha}{1+\alpha}$.

\begin{theorem}
\label{coverage} Under the conditions of Theorem 1 of Marchand and
Strawderman (2006), that is Assumption \ref{assumption}, consider Bayesian credible intervals
$I_{\pi_0, \alpha(\cdot)}$ associated with prior $\pi_0$ and a
spending function $\alpha(\cdot)$ 
such that $\alpha(x)=\alpha$ for all $x$ with $t(x) \leq y_0$.  For the
frequentist coverage $C(\theta) = P_{\theta}(I_{\pi_0,
\alpha(\cdot)}(X) \ni \tau(\theta))$, we then have
\begin{enumerate}
\item[ (a)] $C(\theta) =  \frac{1}{1+\alpha} (> \frac{1-\alpha}{1+\alpha})$ for all $\theta$
such that $\tau(\theta)=0$;
\item[ (b)]  Moreover, we have $C(\theta) >
\frac{1-\alpha}{1+\alpha}$ for all $\theta$ such that
$\tau(\theta) \geq 0$ as long as $\alpha(x)$ satisfies, for all $x$,
\begin{equation}
\label{alphabounds} \frac{(1-\alpha)\,G(-t(x)) +
\frac{\alpha^2}{1+\alpha}}{1-G(-t(x))} \leq \alpha(x) \leq
(\frac{\alpha}{1+\alpha})\frac{1}{1-G(-t(x))} .
\end{equation}
\end{enumerate}
\end{theorem}
{\bf Proof.}
\begin{enumerate}
\item[ {\bf (a)}] First, observe that for $\theta$ such that $\tau(\theta)=0$, the pivot assumption for 
$-T(X,\theta) = \frac{\tau(\theta) - a_1(X)}{a_2(X)}$ implies that $-t(X)=\,-\frac{a_1(X)}{a_2(X)}$ 
has cdf $G$ whenever $\tau(\theta)=0$.   Hence, for $\theta$ such that $\tau(\theta)=0$, we have
$$
P_{\theta}(I_{\pi_0, \alpha(\cdot)}(X) \ni 0)  =
P_{\theta}(\alpha(X)=\alpha)  =  P_{\theta}(t(X) \leq y_0) \\
\;  =  1-G(-y_0)  =  \frac{1}{1+\alpha}\,. $$ 

 \item[ {\bf (b)}]  With the case $\tau(\theta)=0$ addressed in part (a), we consider $\tau(\theta) >0$.  
First, observe that the confidence interval $I_1(X)= [l_1(X), u_1(X)] =
 \max\{0, a_1(X) + a_2(X) G^{-1}(\frac{\alpha}{1+\alpha})\},
 a_1(X) + a_2(X) G^{-1}(\frac{1}{1+\alpha})\}$
has the same frequentist coverage as $I_1^*(X)=[a_1(X) + a_2(X)
G^{-1}(\frac{\alpha}{1+\alpha})\},a_1(X) + a_2(X)
G^{-1}(\frac{1}{1+\alpha})\}]$ equal to
$P_{\theta}(G^{-1}(\frac{\alpha}{1+\alpha}) \leq
\frac{\tau(\theta) - a_1(X)}{a_2(X)} \leq
G^{-1}(\frac{1}{1+\alpha})) = G(G^{-1}(\frac{1}{1+\alpha})) -
G(G^{-1}(\frac{\alpha}{1+\alpha})) = \frac{1-\alpha}{1+\alpha}.$
Now, we show that the given conditions on $\alpha(\cdot)$ imply
that $I_{\pi_0, \alpha(\cdot)} \supseteq I_1$; with the inclusion being strict with probability greater than $0$ for all $\theta$; 
which will lead to the result directly.  Indeed, we have by the upper bound in
(\ref{alphabounds}) and Lemma \ref{atol}: $u_{\alpha(\cdot)}(x)
\geq a_1(x) + a_2(x) G^{-1}(1 - \frac{\alpha}{1+\alpha})=u_1(x)$.
Similarly, from the lower bound (\ref{alphabounds}) and Lemma
\ref{atol} we obtain $l(x) \leq a_1(x) +a_2(x) G^{-1}\{G(-t(x) +
\alpha(1-G(-t(x))) -\frac{\alpha^2}{1+\alpha} - (1-\alpha)
G(-t(x)) \}= a_1(x) + a_2(x)
G^{-1}(\frac{\alpha}{1+\alpha})=l_1(x).$  \qed 
\end{enumerate}

\begin{corollary} 
\label{etcoverage} Under Assumption \ref{assumption}, the $G$-equal-tails credible interval
$I_{\pi_0, \alpha(\cdot)}$, given in Definition \ref{alphaeqt}, has minimum frequentist
coverage $C(\theta)$ greater than $\frac{1-\alpha}{1+\alpha}$ for
all $\theta$ such that $\tau(\theta) \geq 0$.
\end{corollary}
{\bf Proof.}  It suffices to show directly that
(\ref{alphabounds}) is satisfied for the selection
$\alpha(x)=\alpha_{\hbox{eqt}}(x)$ given in (\ref{spendinghpd})
for $x$ such that $t(x) \geq y_0$.  Indeed, we have for such
$x$'s:
$$ \alpha_{\hbox{eqt}}(x) (1-G(-t(x))) = \frac{\alpha}{2} + \frac{1-\alpha}{2} \, G(-t(x))
\leq \frac{\alpha}{2} + \frac{1-\alpha}{2} \, G(-y_0) =
\frac{\alpha}{1+\alpha}\,,$$ and
\begin{eqnarray*}
\alpha_{\hbox{eqt}}(x) (1-G(-t(x))) - (1-\alpha)\, G(-t(x)) -
\frac{\alpha^2}{1+\alpha} &=& \frac{\alpha(1-\alpha)}{2(1+\alpha)}
-
\frac{1-\alpha}{2}\, G(-t(x)) \\
\, &\geq& \frac{\alpha(1-\alpha)}{2(1+\alpha)} -
\frac{1-\alpha}{2}\, G(-y_0))=0.  \;\;\;\;\;\;\;\;\;\;\;\;\; \Box
\end{eqnarray*}

\begin{remark}
In cases where the underlying pivotal distribution is non-symmetric,
Corollary \ref{etcoverage} is a new result, generalizing Theorem 1 of MS (2006), and is widely applicable given the lack of assumptions on $g_0$.  
Also, the bounds of the equal-tails procedure are easier to evaluate than that of the HPD
credible interval.  And the findings of Theorem \ref{coverage} go beyond a single procedure, even in the symmetric case, by providing a class of credible sets, as specified by a spending function, with frequentist coverage bounded below by $\frac{1-\alpha}{1+\alpha}$.  
\end{remark}

We do not have a recommended prescription for the choice of the spending function among those specified by Theorem \ref{coverage} as guaranteeing minimal frequentist coverage of at least $\frac{1-\alpha}{1+\alpha}$.  The $G-$equal-tails choice is simple, intuitively appealing and matches the HPD procedure under symmetry of the pivotal density, while upper tailed and lower tailed choices are not allowed for $x$ such that $t(x) \geq y_0$.   The bounds in (\ref{alphabounds}): (i) $\alpha_1(x) = \frac{(1-\alpha)\,G(-t(x)) + \frac{\alpha^2}{1+\alpha}}{1-G(-t(x))}$   and (ii) $\alpha_2(x) =
(\frac{\alpha}{1+\alpha})\frac{1}{1-G(-t(x))}$, are other interesting choices which push extremally $I_{\pi,\alpha(\cdot)}$ towards
$+\infty$ and $0$ respectively.  Finally, along with these choices, it might be feasible to minimize the length of the credible interval under the restrictions imposed by Theorem \ref{coverage}.

\section{Examples}

At the risk of some redundancy with the examples provided by MS (2006), it is still beneficial here to present various applications with accompanying commentary.  
Assumption 1 is satisfied in all of the examples below with the underlying family of transformations (distributions) being either
the location family, the scale family, or the location-scale family.  In all of the examples, Theorem \ref{coverage} 
and Corollary \ref{etcoverage} provide conditions on the spending function $\alpha(\cdot)$ so that the Bayesian intervals
$I_{\pi_0,\alpha(\cdot)}(X)$ have minimal frequentist coverage greater than $\frac{1-\alpha}{1+\alpha}$ for all $\theta$ 
such that $\tau(\theta) \geq 0$.  These intervals include the equal-tails procedure given in Definition \ref{etdef} and can be evaluated in general using 
the expression given in Lemma \ref{atol}. 

\begin{enumerate}

\item[ (A)] ({\bf a single location parameter}) $X \sim f_0(x-\theta)$; 
$\tau(\theta)=\theta \geq 0$; $T(X,\theta)=X-\theta$;
$\pi_H(\theta)= \mathbb{I}_{\mathbb{R}}(\theta), \pi_0(\theta)=\mathbb{I}_{[0,\infty)}(\theta)$.  In such cases, all Bayes credible sets
$I_{\pi_0,\alpha(\cdot)}$ (with credibility $1-\alpha$), with the spending function $\alpha(\cdot)$ satisfying the conditions of Theorem \ref{coverage} and the bounds in (\ref{alphabounds}), have necessarily minimum frequentist coverage bounded below by $\frac{1-\alpha}{1+\alpha}$.  Through the transformations $X \to X-a$ and $X \to -X+a$, one can reduce all lower bounded restrictions $\theta \geq a$ and upper bounded restrictions $ \theta \leq a$ to the case $\theta \geq 0$ considered here and we will not make further explicit mention of such transformations below.

\end{enumerate}

\begin{remark}
Results such as those in ({\bf A}) are applicable as well for several observations by conditioning on a maximal invariant statistic $V$.  Such a maximal invariant statistic $V$ is an ancillary statistic and specifically an invariant function such that every other invariant statistic is a function of $V$.
Indeed, suppose that $X=(X_1, \ldots, X_n) \sim f_0(x_1-\theta, \ldots, x_n-\theta)$, where $f_0$ is known and where the $X_i$'s are not necessarily independently distributed.  Here, $V=(X_2-X_1, \ldots, X_n-X_1)$ is a maximal invariant statistic.  One can then proceed, for a given value $v$ of $V$, with an interval estimate $I_{\pi_0,\alpha(\cdot,v)}(X_1,v)$ as given in Lemma \ref{atol} with
$G \equiv G_v$ representing the cdf of the pivot $X_1 -\theta$ conditional on $V=v$, and $\alpha(x,v)$ satisfying the conditions of Theorem \ref{coverage} and (\ref{alphabounds}).  This is feasible by the pivot and ancillarity property with the joint distribution of $(X_1-\theta,V)$ independent of $\theta$.  In such a case, Theorem \ref{coverage} applies to the conditional frequentist coverage $C(\theta,v)= P_{\theta}(I_{\pi_0, \alpha(\cdot,v)}(X,v) \ni \tau(\theta)|V=v)$ yielding the inequality $C(\theta,v) > \frac{1-\alpha}{1+\alpha}$ for all $\theta \geq 0$.  Since this is true for all $v$, the unconditional frequentist coverage $C(\theta)$ of the Bayes credible set $I_{\pi_0, \alpha(\cdot, \cdot)}(X,V)$ will also exceed $\frac{1-\alpha}{1+\alpha}$ for all $\theta \geq 0$ (see MS 2006, for more details related to a multivariate Student model).  In the same vein, all the scenarios below ({B to G}), although presented for simplicity in the single observation case, are also applicable in presence of a sample by conditioning on a maximal invariant statistic.
 
\end{remark} 

\begin{enumerate}
\item[ (B)] ({\bf a lower bounded scale parameter}) $X \sim \frac{1}{\theta}
f_1(\frac{x}{\theta}) \, \mathbb{I}_{(0,\infty)}(x)$ with $\theta \geq a$; $\tau(\theta)=\log(\theta) - \log(a) \geq 0$;
$T(X,\theta)= \log(X) - \log (a) - \tau(\theta)$; $\pi_H(\theta)=\frac{1}{\theta}  \mathbb{I}_{(0,\infty)}(\theta)$, $\pi_0(\theta)=
\frac{1}{\theta} \mathbb{I}_{[0,\infty)}(\tau(\theta))$.  Here, an interval estimate of 
$\tau(\theta)$ provides an interval estimate of $\theta$.  Important models include Gamma, Weibull, 
Fisher, among others.  A familiar set-up where the results can be applied arises in random effects analysis of variance
models with a Fisher distributed pivot (see Zhang and Woodroofe, 2002, for details).  As in (A), for a sample $X=(X_1, \ldots, X_n) \sim
\frac{1}{\theta^n} \, f_1(\frac{x_1}{\theta}, \ldots, \frac{x_n}{\theta}) \prod_i \mathbb{I}_{(0,\infty)}(x_i)$, Theorem \ref{coverage} and Corollary \ref{etcoverage} are applicable by conditioning
on the maximal invariant statistic $V=(\frac{X_1}{X_n}, \ldots, \frac{X_{n-1}}{X_n})$. 

\end{enumerate}

\begin{remark}
Further applications consist of power parameter families where we have a scale family for an observable $Y$ and the model of interest are the distributions for $X=e^Y$.  As as simple illustration, consider the Pareto model for $X$ with densities $ \frac{\gamma}{x^{\gamma+1}} 1_{(1,\infty)}(x)$ and the parametric constraint $\gamma \in (1,\gamma_0)$.  In such cases, we have
that $\gamma_0 \log(X) \sim \hbox{Exp}(\theta)$ with $\theta=\frac{\gamma_0}{\gamma} \geq 1$ and the results in (B) apply.
 
\end{remark} 

\begin{enumerate}
\item[ (C)] ({\bf location-scale families }) $(X_1,X_2) \sim
\frac{1}{\theta_2^2} \, f_2(\frac{x_1-\theta_1}{\theta_2},\frac{x_2}{\theta_2}) \, \mathbb{I}_{(0,\infty)}(x_2)$;
$\tau(\theta)=\theta_1 \geq 0$; $T(X,\theta)={X_1-\theta_1 \over
X_2}$; $\pi_H(\theta)=\frac{1}{\theta_2} \mathbb{I}_{(0,\infty)}(\theta_2)
\mathbb{I}_{(-\infty,\infty)}(\theta_1)$, $\pi_0(\theta)=\frac{1}{\theta_2} \mathbb{I}_{(0,\infty)}(\theta_2)
\mathbb{I}_{[0,\infty)}(\theta_1)$.  This set-up encompasses, but is not limited to, the basic normal 
case: $Y_1, \ldots Y_n \sim^{ind.} N(\mu, \sigma^2)$ with $\sigma^2$ unknown and $\mu \geq 0$, and by taking $X_1$ and $X_2$ 
respectively as the sample mean and standard deviation of the $Y_i$'s.  More generally, the results apply
for linear models $Y = Z\beta + \epsilon$,
$\epsilon \sim N(0,\sigma^2 I_n)$ where the objective is to estimate a lower-bounded linear combination
$\tau(\theta)=l'\beta$,  by setting $X_1=\hat{\beta}(Z'Z)^{-1}Z'Y$, $X_2^2=\|Y-Z\beta\|^2$,
$\theta_1=\beta$, $\theta_2=\sigma$.  Here, the pivot $T(X,\theta)$ has a Student distribution.
Alternatively, if the objective is to estimate a lower bounded scale $\theta_2$, one can proceed as in (B).

\item[ (D)] ({\bf linear combination of several location parameters}) $X=(X_1, \ldots, X_p)
\sim f_0(x_1-\theta_1, \ldots, x_p-\theta_p)$; $\tau(\theta)=
\sum_{i=1}^p a_i \theta_i$; $\pi_H(\theta)= \mathbb{I}_{\mathbb{R}^p}(\theta), \pi_0(\theta)=\mathbb{I}_{[0,\infty)}(\tau(\theta))$,
$T(X,\theta)= (\sum_{i=1}^p a_i X_i) -
\tau(\theta)$.  This set-up includes, for instance, estimating a difference $\theta_1 - \theta_2$ with an order constraint
$\theta_1 \geq \theta_2$.

\item[(E)] ({\bf multivariate location-scale families with homogeneous scale} ) \\
In (D), we can incorporate a common scale and apply the results of this paper for estimating 
a lower bounded linear combination with    
$X=(X_1, \ldots, X_p\,,\,X_{p+1}) \sim
f_0(\frac{x_1-\theta_1}{\theta_{p+1}}, \ldots, \frac
{x_p-\theta_p}{\theta_{p+1}}, \frac{x_{p+1}}{\theta_{p+1}})$, 
$\tau(\theta)= \sum_{i=1}^p a_i \theta_i$, $T(X,\theta)= {
(\sum_{i=1}^p a_iX_i) - \tau(\theta) \over X_{p+1}}$, and 
$\pi_0(\theta)=\frac{1}{\theta_{p+1}}1_{(0,\infty)}(\theta_{p+1})
1_{[0,\infty)}(\tau(\theta))$.  

\item[ (F)] ({\bf several scale parameters }) \\
$(X_1, \ldots, X_p) \sim
(\Pi_{i=1}^{p} \frac{1}{\theta_i})\; f_1(\frac{x_1}{\theta_1},
\ldots, \frac{x_p}{\theta_p})$; $\tau(\theta)= \sum_{i=1}^p a_i
\log (\theta_i)$, $\pi_H(\theta)=\prod_i \frac{1}{\theta_i}
\mathbb{I}_{(0,\infty)} (\theta_i)$, $\pi_0(\theta)= \pi_H(\theta) \, \mathbb{I}_{[0,\infty)}(\tau(\theta))$.  This can consist, for instance with $p=2$, $a_1=1, a_2=-1$, of estimating a lower bounded ratio $\frac{\theta_2}{\theta_1} \geq 1$ 
of two scale parameters.

\item[ (G)] ({\bf quantiles in location-scale families}) $X_i \sim^{\hbox{ind.}}
N(\mu,\sigma)$, $ i=1,\ldots n$, $\theta=(\mu,\sigma)$, $\tau(\theta)= \mu + \eta \sigma \geq 0$, $\pi_H(\theta)=\frac{1}{\theta_2} \mathbb{I}_{(0,\infty)}(\theta_2)
\mathbb{I}_{(-\infty,\infty)}(\theta_1)$,
$\pi_0(\mu, \sigma)=\frac{1}{\sigma} \mathbb{I}_{(0,\infty)}(\sigma)
\mathbb{I}_{[0,\infty)} (\mu + \eta \sigma)$.  $T(X,\theta)=
\frac{\bar{X}-\mu -\eta\sigma}{S}$.  Here, $T(X,\theta)$ is distributed as non-central Student. 
The applications are not restricted to normality and are applicable in general for location-scale families as in (C).
\end{enumerate}

\section{Concluding remarks}

For a large variety of situations with a lower bounded parametric constraint, we have obtained a class of Bayesian $(1-\alpha) \times 100\%$ 
credible sets which provide minimal frequentist probability coverage exceeding $\frac{1-\alpha}{1+\alpha}$.  These Bayesian confidence intervals include an equal tailed modification or approximation of the HPD credible set which coincides with the latter when the distribution of the underlying pivot is symmetric.   In non-symmetric cases not covered by Marchand and Strawderman (2006), our findings provide instances of Bayesian credible sets with given minimal frequentist coverage and hence fill a gap in the literature.   In comparison to earlier results for normal models, as well as the symmetric models considered by Marchand and Strawderman (2006), the findings here relative to the HPD are not new, but those related to other Bayesian credible sets are an addition.
In seeking to evaluate the frequentist performance of Bayesian confidence intervals, our results illustrate that the choice of bounds or spending function matters, so that there does not necessarily exist a single universal assessment of their frequentist performance even in a given specific problem.  

\section*{Acknowledgements}
The authors are grateful to two reviewers, an associate editor, and the editor for useful comments and suggestions which led to a more self-contained and readable
manuscript.   Eric Marchand's research is supported in part by a grant from the Natural Sciences and Engineering Research Council of Canada, and William Strawderman's research is partially supported by a grant from the Simons Foundation (\#209035).

\renewcommand{\baselinestretch}{1.2}
\section*{References}
\small
\medskip\noindent
Berger, J.O. (1985).  {\it Statistical Decision Theory and
Bayesian Analysis}.  Springer-Verlag, New York, 2nd edition.
\medskip\noindent

\medskip\noindent
Eaton, M. L. (1989).  Group invariance Applications in Statistics. NSF-CBMS Regional Conference
Series in Probability and Statistics 1. Hayward, CA.

\medskip\noindent
Feldman, G.J. and Cousins, R. (1998).  Unified approach to the
classical statistical analysis of small signals.  {\it Physical
Review D}, {\bf 57}, 3873-3889.

\medskip\noindent
Fraser, D.A.S. (2011).  Is Bayes posterior just quick and dirty confidence?
{\it Statistical Science}, {\bf 26}, 299-316.

\medskip\noindent
Mandelkern, M. (2002).  Setting Confidence Intervals for Bounded
Parameters with discussion.  {\it Statistical Science}, {\bf 17},
149-172.

\medskip\noindent
Marchand, \'{E}., Strawderman, W. E., Bosa, K., and Lmoudden, A.
(2008).  On the frequentist coverage of Bayesian credible
intervals for lower bounded means. {\it Electronic Journal of
Statistics}, {\bf 2}, 1028-1042.

\medskip\noindent
Marchand, \'{E}. and Strawderman, W. E. (2006). On the behaviour
of Bayesian credible intervals for some restricted parameter space
problems. {\it Recent Developments in Nonparametric Inference and
Probability : A Festschrift for Micheal Woodroofe,}, IMS Lecture
Notes-Monograph Series, {\bf 50}, pp. 112-126.


\medskip\noindent
Robert,  C.P. (2011).  Discussion of ``Is Bayes posterior just quick and dirty confidence?''
by D.A.S. Fraser. {\it Statistical Science}, {\bf 26}, 317-318.

\medskip\noindent
Roe, B. and Woodroofe, M. (2000). Setting confidence belts. {\it
Physical Review D}, {\bf 63}, 013009/01-09.

\medskip\noindent
Wasserman, L.  (2011).  Frasian inference. {\it Statistical Science}, {\bf 26}, 322-325.

\medskip\noindent
Zhang, T. and Woodroofe, M. (2003).  Credible and confidence sets
for restricted parameter spaces. {\it Journal of Statistical
Planning and Inference}, {\bf 115}, 479-490.

\medskip\noindent
Zhang, T. and Woodroofe, M. (2002).  Credible and confidence sets
for the ratio of variance components in the balanced one-way
model. {\it Sankhy$\bar{a}$: Special issue in memory of D. Basu},
{\bf 64}, 545-560.

\end{document}